






\documentstyle{amsppt}
\magnification=\magstep1
\NoRunningHeads
\NoBlackBoxes
\def\Ker{\text{Ker }}
\def\dim{\text{dim }}

\def\Coker{\text{Coker }}
\def\Ker{\text{Ker }}
\vsize=7.4in

\topmatter


\font\title=cmr10 scaled\magstep1
\centerline{\title STABILITY RESULTS ON INTERPOLATION SCALES}
\vskip 0.15in
\centerline{\title OF QUASI-BANACH SPACES AND APPLICATIONS}

\leftheadtext{Nigel Kalton and Marius Mitrea}
\rightheadtext{Stability results on interpolation scales}

\author{Nigel Kalton and Marius Mitrea}\endauthor

\thanks 1991 {\it Mathematics Subject Classification.} Primary 46A16, 42B20,
Secondary 47A53, 46B70.
\endthanks

\abstract{We investigate stability of Fredholm properties on interpolation
 scales
of quasi-Banach spaces. This analysis is motivated by problems arising in PDE's
and several applications are presented.}\endabstract

\endtopmatter

\document

\heading\S 1.\quad Introduction
\endheading

In this paper we initiate the study of stability of Fredholm properties of
operators on complex interpolation scales of quasi-Banach spaces. As such, this
is a natural continuation and extension of
previous work in the literature (cf., e.g., the articles
\cite{32}, \cite{35}, \cite{7}, \cite{8})
which only deals with the case of Banach spaces.

The first task is to identify the essential functional analytic elements such
 that
a satisfactory stability theory can be developed in the context of quasi-Banach
 spaces.
One of the main problems is that we are forced to work with non-locally convex
spaces and, consequently, several central theorems from the Banach space
setting
are no longer available (most notably duality results based on the Hahn-Banach
 theorem).
Let us point out that all existing papers addressing similar issues
employ in one form or another results from Banach space theory which do not
 carry over
to the quasi-Banach context (for instance, it is assumed that the dual of the
interpolation scale is an interpolation scale itself).
Thus, it is desirable to device a different approach which, in particular,
 avoids
the use of duality.

In section 2 we indicate that this can be accomplished in a rather general,
flexible framework. For instance, here we show that for an interpolation
operator being an isomorphic embedding
or being onto are stable properties on interpolation scales of quasi-Banach
 spaces
(cf. Theorem 2.3, Theorem 2.5). The same applies to the property of being
Fredholm and to the index of the operator (see Theorem 2.9).
In fact the same results can be extended to analytic families of operators
(compare with \cite{11}). Let us also note that our
approach is constructive in the sense that all perturbation constants
explicitly
depend on the original ones (this is relevant in a number of applications).

Section 3 contains a discussion of the complex method applied to
interpolation couples of quasi-Banach spaces. Among other things,
here we examine the role played by the Maximum Modulus Principle in the context
of analytic functions with values in quasi-Banach spaces.

Our primary motivation for studying scales of quasi-Banach spaces
steams out from certain PDE's problems which involve the family
of (real variable) $H^p$ spaces for $p\leq 1$. The point is that many
 perturbations
of, say, $H^1$-invertibility results are of considerable interest.

In order to make this idea more transparent, recall
first an old result due to S. Zaremba (\cite{36}) to the effect that the unique
solution of the Dirichlet problem for the Laplacian in Lipschitz domains with
 H\"older
continuous data extends to a H\"older continuous function on the closure of the
 domain
(with the same, small exponent). The original proof is based on the
construction
of barriers and now the question arises whether
it is possible to represent the solution in the layer potential form
with a H\"older continuous boundary density (much in the spirit of the $L^2$
 theory
in \cite{34}).

In proving this, a key step is to show that $\frac{1}{2}I+{\Cal K}^*$, the
 adjoint of the
double layer potential operator, is an isomorphism of $H^p$ for
$1-\epsilon<p<1$
 for
some small $\epsilon >0$.
This has been first proved by R. Brown \cite{4} (at least for starlike
 Lipschitz
domains) by a careful analysis of the Neumann
and Regularity problems for the Laplacian with data in $H^p$, $1-\epsilon<p<1$.
Nonetheless, it is natural to investigate if this result can be obtained more
 directly
by combining a theorem of B. Dahlberg and C. Kenig (\cite{14}) concerning the
invertibility of $\frac{1}{2}I+{\Cal K}^*$ on $H^1$ with a general stability
result on the scale $\{H^p\}_{0<p<\infty}$.
In fact, this particular question has also been posed to us by Eugene Fabes.
To this, we give a positive answer in the final section of the paper.

Our method of proof is quite flexible and
here we also show that similar considerations apply to other situations as
well.
In particular, this is the case for the three
 dimensional
system of elastostatics in Lipschitz domains, for which we are able to further
 refine
the theory developed in \cite{15}, \cite{16}.

\vskip 0.08in
\noindent{\it Acknowledgments.} The second named author would like to thank
Eugene Fabes for several discussions and his interest in this work.

\heading\S 2.\quad Interpolation scales of quasi-Banach spaces. The abstract
 setup
\endheading

Throughout the paper, for a quasi-normed space $(X,\|\cdot\|_X)$,
we shall denote by $\rho=\rho(X)$ its {\it modulus of concavity},
i.e. the smallest positive constant for which
$$
\|x+y\|_X\leq\rho(\|x\|_X+\|y\|_X),\qquad x,y\in X
$$
(note that always $\rho(X)\geq 1$).  We recall the Aoki-Rolewicz theorem
\cite{26},\cite{31}, which asserts that $X$ can be given an equivalent
$r$-norm (where $2^{1/r-1}=\rho$) i.e. a quasi-norm which also satisfies
the inequality:
$$ \|x+y\| \leq (\|x\|^r+\|y\|^r)^{1/r}.$$
In general a quasi-norm need not be continuous but an $r$-norm is
continuous. We shall assume however, throughout the paper that all
quasi-norms considered are continuous: in fact, of course it would
suffice to consider an $r$-norm for suitable $r$.

Let $U$ be a fixed, (Hausdorff) locally compact, locally connected topological
 space,
referred to in the sequel as the {\it space of parameters}, and let
$\delta :U\times U\to\bold C$ be a (fixed) continuous function such
that $\delta(z,w)=0$ if and only if $z=w$. Also,
suppose that $Z$ is a complex (Hausdorff) topological vector space, called the
{\it ambient space}.

A family $\Cal F$ of functions
which map $U$ into $Z$ is called {\it admissible} (relative to $U$, $\delta$
and $Z$) provided the following axioms are satisfied:
\roster
\item ${\Cal F}$ is a (complex) vector space endowed with a quasi-norm
$\|\cdot\|_{\Cal F}$ with respect to which this is complete (i.e. ${\Cal F}$ is
a quasi-Banach space);
\item the point-evaluation mappings $\text{ev}_w:{\Cal F}\to Z$, $w\in U$,
defined by $\text{ev}_w(f):=f(w)$ are continuous;
\item for any $K$ compact subset of $U$ there exists a positive constant
$D_K({\Cal F}$
such that for any $w\in K$ and any $f\in{\Cal F}$ with $f(w)=0$, it then
follows
that the mapping $U\setminus\{w\}\ni z\mapsto f(z)/\delta(z,w)\in Z$
extends to an element in ${\Cal F}$ and
$$
\left\|\frac{f(z)}{\delta(z,w)}\right\|_{\Cal F}\leq D_K({\Cal F})\|f\|_{\Cal
F}.
\tag{2.1}
$$
\endroster
We remark that, occasionally, we will write $D_K$ for $D_K({\Cal F})$,
and $D_w$ for $D_{\{w\}}$.

Assume that such an admissible family ${\Cal F}$
has been fixed (we shall also call the elements of ${\Cal F}$ admissible
 functions).
Introduce the {\it intermediate spaces} $X_w$, $w\in U$, by
$$
X_w:=\{x\in Z;\,x=f(w)\,\,\text{for some}\,\,f\in{\Cal F}\}\subseteq Z.
\tag{2.2}
$$

Under these hypotheses it is trivial to establish:

\proclaim{Lemma 2.1} For each $w\in U$, the
intermediate space $X_w$ becomes a quasi-Banach space when equipped with
$\|x\|_w:=\text{inf}\,\{\|f\|_{\Cal F};\,f\in{\Cal F},\,f(w)=x\}$, $x\in X_w$.

Furthermore, the modulus of concavity of the interpolation quasi-norm
$\|\cdot\|_w$ does not exceed $\rho({\Cal F})$,
and the inclusion $X_w\hookrightarrow Z$ is continuous.
\endproclaim

\demo{Proof}We need only remark that $X_w$ is naturally identified with
the quotient space $\Cal F/\Ker (ev_{w})$.
\hfill$\blacksquare$
\enddemo

The collection of all intermediate spaces $\{X_w\}_{w\in U}$ is called {\it an
interpolation scale} (relative to the space of parameters $U$, the function
$\delta$, the ambient space $Z$, and the admissible family
${\Cal F}$; however, occasionally we may refer to $(U,\delta,Z,{\Cal F})$ as
being an interpolation scale).

Consider next two interpolation scales: $\{X_w\}_{w\in U}$,
relative to some ambient space $Z$ and some admissible family ${\Cal F}$, and
$\{Y_w\}_{w\in U}$,
relative to another ambient space $W$ and some admissible family ${\Cal G}$
(note that $U$ and $\delta$ are assumed to be same for both scales).
A family of operators $\{T_w\}_{w\in U}$,
$$
T_w:X_w\to Z,\qquad w\in U,
$$
will be called {\it an interpolating  family of operators},
provided the application
$$
{\Cal F}\ni f\mapsto [U\ni w\mapsto T_w(f(w))]\in{\Cal G}
$$
is well-defined, linear and satisfies
$$
\|T_{.}(f(\cdot))\|_{\Cal G}\leq M\,\|f\|_{\Cal F},\quad f\in{\Cal F},
\tag{2.3}
$$
for some positive, finite constant $M$. The best constant in (2.3) is called
the
{\it interpolation norm} of $\{T_w\}_{w\in U}$.

Note that, if a single operator
$$
T:\bigcup_{w\in U}X_w\to\bigcup_{w\in U}Y_w
$$
maps admissible functions into admissible functions linearly and boundedly,
then the
``constant'' family $T_w:=T|_{X_w}$, $w\in U$, becomes an interpolation family
of
analytic operators. In the sequel, we shall call such a mapping $T$
{\it an interpolation operator}.

\proclaim{Lemma 2.2} Assume that $\{T_w\}_{w\in U}$ is an interpolating  family
of operators between two
interpolation scales $\{X_w\}_{w\in U}$ and $\{Y_w\}_{w\in U}$ as before.
Then, for each $w\in U$, the operator $T_w$ maps the intermediate space $X_w$
linearly and boundedly into $Y_w$, with norm not exceeding the
interpolation norm of the family $\{T_w\}_{w\in U}$.
\endproclaim

\demo{Proof} The argument is straightforward. Let $M$ estimate
the interpolation norm of $\{T_w\}_w$
and, for $x\in X_w$, $\epsilon >0$, pick $f\in{\Cal F}$ such that $f(w)=x$ and
$\|f\|_{\Cal F}\leq\|x\|_w+\epsilon$. Then, since $T_{.}(f(\cdot))\in{\Cal G}$,
it
follows that $T_wx\in Y_w$ and
$$
\|T_wx\|_w\leq\|T_{.}(f(\cdot))\|_{\Cal G}\leq M\|f\|_{\Cal F}\leq
M(\|x\|_w+\epsilon).
$$
The conclusion follows.
\hfill$\blacksquare$
\enddemo

To state the main results of this section, we still need to introduce some more
notation. Consider $(X,\|\cdot\|_X)$,
$(Y,\|\cdot\|_Y)$ two quasi-Banach spaces and let $T:X\to Y$ be a linear,
 bounded
operator.
Define $\kappa(T;\,X,Y)$ to be the smallest constant so that if $y\in Y$ then
 there
exists $x\in X$ so that $Tx=y$ and $\|x\|_X\le \kappa(T;\,X,Y)\|y\|_Y$.
Note that, by the Open Mapping Theorem (which remains valid in the context of
quasi-Banach spaces; cf. Theorem 1.4 in \cite{26}), $\kappa(T;\,X,Y)$ is finite
if and only if $T$ maps $X$ onto $Y$.

We also let $\eta(T;\,X,Y)$ be the largest constant so that
$\eta(T;\,X,Y)\|x\|_X\leq \|Tx\|_Y$ for each $x\in X$. Once again by virtue
of the Open mapping Theorem, $\eta(T;\,X,Y)$ is strictly positive if and only
if
the operator $T$ is injective and has closed range.

Finally, if $\{T_w\}_{w\in U}$ is an interpolating analytic family of operators
acting between two interpolation scales $\{X_w\}_{w\in U}$ and $\{Y_w\}_{w\in
U}$
then, for each $w\in U$, we abbreviate $\kappa(T_w;\,X_w,Y_w)$
and $\eta(T_w;\,X_w,Y_w)$ by $\kappa(T_w)$ and by $\eta(T_w)$, respectively.

Next we can state our first stability result, concerning the ontoness of an
interpolation operator.

\proclaim{Theorem 2.3} Assume that $\{T_w\}_{w\in U}$ is an interpolation
analytic family of operators between two interpolation scales $\{X_w\}_{w\in
U}$
and $\{Y_w\}_{w\in U}$.
Furthermore, assume that $\kappa(T_{w_0})$ is finite for some ${w_0}\in U$.
Then there exist a neighborhood $\omega$ of $w_0$ in $U$ and a finite
positive constant $C$ such that $\kappa(T_z)\leq C$ for each $z\in\omega$.
\endproclaim

In order to prove this theorem we shall need a lemma, which is in fact a
part of the standard proof of the Open Mapping Theorem.

\proclaim{Lemma 2.4} Suppose that $(X,\|\cdot\|_X)$,
$(Y,\|\cdot\|_Y)$ are two quasi-normed spaces such that $X$ is complete.
Also, suppose that $T:X\to Y$ is a linear, bounded operator for which the
following property is true: there exist $0<C_0<+\infty$ and $0<a<1$
such that for each $y$ in the unit sphere of $Y$ one can find $x\in X$ with
$\|x\|_X\leq C_0$ and $\|y-Tx\|_Y\leq a$.

Then $T$ is onto and $\kappa(T;\,X,Y)\leq C_1$ for some $C_1$ depending
exclusively on $C_0$, $\rho(X)$ and $a$.
\endproclaim

\demo{Proof} For fixed $y\in Y$, it one can inductively define
$(x_n)_{n=1}^{\infty}$ in $X$ so that $\|y-\sum_{k=1}^nTx_k\|_Y \le
a^n\|y\|_Y$ and $\|x_n\|_X\le C_0a^n\|y\|_Y$.  Then if
$2^{1/r-1}=\rho(X)$ we observe that
$(\sum_{n=1}^{\infty}\|x_n\|_X^r)^{1/r}\le C_0(1-a^r)^{-1/r}\|y\|_Y$.
Hence $\sum_{n=1}^{\infty}x_n=x$ converges and $\|x\|_X\le C_1\|y\|_Y$
where $C_1$ depends on $C_0,\rho(X)$ and $a$.  Clearly $Tx=y$.
\hfill$\blacksquare$
\enddemo

We are now in a position to present the proof of Theorem 2.3.

\demo{Proof of Theorem 2.3} Let ${\Cal F}$, ${\Cal G}$ be the two classes of
admissible functions corresponding to the two interpolation scales
$\{X_w\}_{w\in U}$
and $\{Y_w\}_{w\in U}$.
Also, denote by $\rho({\Cal F})$, $\rho({\Cal G})$, respectively, their moduli
of
concavity, and by $M$ the interpolation norm of the family $\{T_w\}_w$.

Now, choose a neighborhood $\omega\subseteq U$ of $w_0$ such that
$|\delta(z,w_0)|$ is sufficiently small, uniformly for $z\in\omega$
(this is made precise shortly),
and let $y$ be an arbitrary vector in $Y_z$ with $\|y\|_z=1$.
Pick $g\in{\Cal G}$ such that $g(z)=y$,
$\|g\|_{\Cal G}\leq 2$, and set $\tilde{y}:=g(w_0)\in Y_{w_0}$. Note that
$\|\tilde{y}\|_{w_0}\leq\|g\|_{\Cal G}\leq 2$. Then, by hypothesis, there
exists
$\tilde{x}\in X_{w_0}$ such that $T_{w_0}\tilde{x}=\tilde{y}$ and
$\|\tilde{x}\|_{w_0}\leq\kappa(T_{w_0})\|\tilde{y}\|_{w_0}$.
Thus, $\|\tilde{x}\|_{w_0}\leq 2\kappa(T_{w_0})$ and,
therefore, there exists $f\in{\Cal F}$
so that $f(w_0)=\tilde{x}$ and $\|f\|_{\Cal F}\leq 2\kappa(T_{w_0})+1$.

Define now $h:U\to W$ (recall that $W$ is the ambient space for the
interpolation scale $\{Y_w\}_{w\in U}$) by
$$
h(w):=\frac{T_wf(w)-g(w)}{\delta(w,w_0)}.
$$
By axiom (3), $h$ belongs to ${\Cal G}$ and
$$
\align
\|h\|_{\Cal G}\leq &D_{w_0}\rho({\Cal G})(\|T_{.}(f(\cdot))\|_{\Cal
G}+\|g\|_{\Cal
G})\leq D_{w_0}\rho({\Cal G})(M\|f\|_{\Cal F}+\|g\|_{\Cal G})\\
\leq &D_{w_0}\rho({\Cal G})(2M\kappa(T_{w_0})+M+2)
=:C(M,w_0,\rho({\Cal G}),\kappa(T_{w_0}))<+\infty .
\endalign
$$
The idea is now that, if $x:=f(z)\in X_z$, then
$\|x\|_z\leq\|f\|_{\Cal F}\leq 2\kappa(T_{w_0})+1$ and
$$
\|T_zx-y\|_z\leq|\delta(z,w_0)|\|h(z)\|_z\leq |\delta(z,w_0)|\|h\|_{\Cal G}
\leq |\delta(z,w_0)|C(M,w_0,\rho({\Cal G}),\kappa(T_{w_0})).
$$
Consequently, if $\omega$ is chosen so that for each $z\in\omega$ one has
$$
|\delta(z,w_0)|C(M,w_0,\rho({\Cal G}),\kappa(T_{w_0}))\le\tsize{\frac{1}{2}},
$$
Lemma 2.4 can be used to conclude the proof of the theorem.
\hfill$\blacksquare$
\enddemo

Our next result deals with the stability of the property of being
an isomorphic embedding for an interpolation operator.

\proclaim{Theorem 2.5} Let $\{T_w\}_{w\in U}$ be an interpolating  family
of operators between two interpolation scales $\{X_w\}_{w\in U}$
and $\{Y_w\}_{w\in U}$, and
assume that $\eta(T_{w_0})$ is strictly positive at some point ${w_0}\in U$.
Then there exist a neighborhood $\omega$ of $w_0$ in $U$ and a
constant $C$ such that $0<C\leq\eta(T_z)$ for each $z\in\omega$.
\endproclaim

First we shall need a technical lemma, essentially to the effect that if an
 admissible
function is large at a certain point in the space of parameters then it remains
large when evaluated at neighboring points.

\proclaim{Lemma 2.6} Consider $\{X_w\}_{w\in U}$ an interpolation scale
(relative to some $\delta$, $Z$, ${\Cal F}$). Let $K$ be a compact subset
of $U$ and suppose $w_0\in U$. Then for any $f\in\Cal F$ and $z_0\in K$
such that $\rho({\Cal F})D_K({\Cal F})|\delta(w_0,z_0)|\le 1$ we have
$$
\|f(z_0)\|_{z_0} \ge \frac1{3\rho({\Cal F})}\|f(w_0)\|_{w_0}-|\delta(w_0,z_0)|
\frac{\rho({\Cal F})}{3}D_K({\Cal F})\|f\|_{\Cal F}.
$$
\endproclaim

\demo{Proof}First pick any $g\in\Cal F$ so that $\|g\|_{\Cal F}\le
2\|f(z_0)\|_{z_0}$ and $g(z_0)=f(z_0)$. Then we may find $h\in\Cal F$ so
that
$\delta(z,z_0)h(z)=f(z)-g(z)$ and $\|h\|_{\Cal F}\le D_K\|f-g\|_{\Cal F}
\le D_K\rho(\|f\|_{\Cal F} +2\|f(z_0)\|_{z_0})$ where $\rho=\rho(\Cal F)$.
It follows that
$$
\|f(w_0)-g(w_0)\|_{w_0}\le |\delta(w_0,z_0)|\|h\|_{\Cal F}
$$
and, hence,
$$
\|g(w_0)\|_{w_0}\ge
\rho^{-1}\|f(w_0)\|_{w_0}-|\delta(w_0,z_0)|\|h\|_{\Cal F}.
$$
Thus,
$$
2\|f(z_0)\|_{z_0} \ge \rho^{-1}\|f(w_0)\|_{w_0}
-|\delta(w_0,z_0)|D_K\rho(\|f\|_{\Cal F}+2\|f(z_0)\|_{z_0}).
$$
Hence, if $D_K\rho|\delta(w_0,z_0)|\le 1$ we have
$$
\|f(z_0)\|_{z_0} \ge \frac1{3\rho}\|f(w_0)\|_{w_0}-|\delta(w_0,z_0)|
\frac{D_K\rho}{3}\|f\|_{\Cal F}.
$$
This completes the proof of the lemma.
\hfill$\blacksquare$
\enddemo

Next we are ready to present the proof of Theorem 2.5.

\demo{Proof of Theorem 2.5}
Let $K$ be a fixed compact neighborhood of $w_0$. Suppose $z_0\in K$
is such that $D_K(\Cal G)\rho(\Cal G)|\delta(w_0,z_0)|\le 1$.
Now pick $x\in X_{z_0}$ with $\|x\|_{z_0}=1$.
We can pick $f\in\Cal F$ so that $\|f\|_{\Cal F}\le 2$
and $f(z_0)=x$. If we set $g(z):=T_z(f(z))$ then $g\in\Cal G$ and, if $M$
stands for
the interpolation norm of the family $\{T_w\}_w$, then $\|g\|_{\Cal G}\le 2M$.
Hence, by Lemma 2.6,
$$
\|T_{z_0}x\|_{z_0} \ge \frac{1}{3\rho(\Cal G)}\|T_{w_0}(f(w_0))\|_{w_0}
-2M\frac{D_K(\Cal G)\rho(\Cal G)}{3}|\delta(w_0,z_0)|.
$$
It remains to estimate $\|T_{w_0}f(w_0)\|_{w_0}$. We have (again by Lemma 2.6)
$$
\align
 \|T_{w_0}f(w_0)\|_{w_0}&\ge \eta(T_{w_0})\|f(w_0)\|_{w_0}\\
&\ge \eta(T_{w_0})\left(\frac{1}{3\rho(\Cal F)}
-|\delta(z_0,w_0)|\frac{D_K(\Cal F)\rho(\Cal F)}{3}\right).
\endalign
$$
Substituting back it is clear that we find a neighborhood $\omega$ or
$w_0$ so that $\eta(T_{z_0})>\frac1{10}(\rho(\Cal F)\rho(\Cal
G))^{-1}\eta(T_{w_0})$
for $z_0\in \omega$.
\hfill$\blacksquare$
\enddemo

\demo{Remark} Essentially the same proof above shows that operator
inequalities of the form $\|x\|_w\leq C'\|S_w x\|_w\leq C''\|T_w x\|_w$,
$x\in X_w$, for $\{S_w\}_w$, $\{T_w\}_w$ two interpolating
families of operators on interpolation scales of quasi-Banach spaces
are stable (in the natural, obvious sense).
In fact, if $\{S_w\}_w$, $\{T_w\}_w$ are two interpolating  families of
operators on interpolation scales of quasi-Banach
spaces such that, for some $w_0$,
$\|S_{w_0}x\|_{w_0}\leq C\|T_{w_0}x\|_{w_0}$, for all $x\in X_{w_0}$,
then for any $\epsilon>0$ there exist $C'>0$ and a neighborhood $\omega$ of
$w_0$ such that
$$
\|S_{w}x\|_w\leq C'\|T_{w}x\|_w+\epsilon\|x\|_w
$$
for any $w\in\omega$ and all $x\in X_w$.
Applications of this type of result to PDE's can be found in, e.g., \cite{30}.
\hfill$\square$
\enddemo

An immediate corollary of Theorem 2.3 and Theorem 2.5 is the stability on
interpolation scales of the property of being an isomorphism. Hence, for
further reference we note the following.

\proclaim{Theorem 2.7} Assume that $\{T_w\}_{w\in U}$ is an interpolating
family of operators between two interpolation scales $\{X_w\}_{w\in U}$
and $\{Y_w\}_{w\in U}$. Then the set of points $z\in U$ at which $T_z:X_z\to
Y_z$ is
an isomorphism is open in $U$.
\endproclaim

Our next goal is to further extend this theorem by studying the stability of
Fredholm properties on interpolation scales. We debut with a preliminary result
in this direction.

\proclaim{Lemma 2.8} Consider $\{T_w\}_{w\in U}$ an interpolating  family
of operators between two interpolation scales $\{X_w\}_{w\in U}$
and $\{Y_w\}_{w\in U}$, and assume that $\kappa(T_z)$ is bounded uniformly for
$z$ in some connected compact subset $K$ of $U$. Then the dimension of
$\text{Ker}\,T_z$ is constant on $K$ (eventually this can be $+\infty$).
\endproclaim

\demo{Proof} The idea of proof is to show that
for each fixed $w_0\in K$ there exists some positive $\epsilon$ such that if
$n\in \bold N$ has $\text{dim}\,\text{Ker}\,T_{w_0}\geq n$, then
also $\text{dim}\,\text{Ker}\,T_z\geq n$ for any $z\in K$ in a small
neighborhood
of $w_0$. Then the conclusion follows from the usual argument.

To this end, assume that $\kappa(T_z)$, $z\in K$, is bounded above by some
finite,
positive $\kappa$, and recall the positive constant $D_K$ corresponding to $K$
as in
axiom (3). Also, fix $w_0\in K$, $n\in\bold N$ with
$\text{dim}\,\text{Ker}\,T_{w_0}\geq n$,
and then pick an $n-$dimensional subspace of $\text{Ker}\,T_{w_0}$.
We shall consider next an Auerbach basis of this subspace, i.e.
$x_1,\ldots,x_{n}\in\text{Ker}\,T_{w_0}\subseteq X_{w_0}$ so that
$\|x_i\|_{w_0}=1$ and $\|\sum a_ix_i\|_{w_0}\ge \max|a_i|$. This idea
has been first employed by Cao-Sagher \cite{8} in the context of Banach
spaces. The existence of Auerbach bases for real quasi-normed spaces
is obtained by maximizing the determinant just as in the real normed case
and will work with only the assumption that the unit ball is closed (e.g. if
the quasi-norm is continuous as we assume). The only modification
necessary for complex scalars is that one seeks the point where the determinant
is
real and has maximum absolute value. We refer to \cite{28}, p.16; note
that the argument in \cite{31}, p. 299 for the complex case is apparently
incorrect.

Pick $f_i\in\Cal F$ with $f_i(w_0)=x_i$ and $\|f_i\|_{\Cal F}\leq 2$,
$1\leq i\leq n$. Then, for each $i$, we have
$U\ni z\mapsto(T_zf_i(z))\delta(z,w_0)^{-1}\in{\Cal G}$ and
$\|(T_{.}(f_i(\cdot)))\delta(\cdot,w_0)^{-1}\|_{\Cal G}\leq 2MD_K$, where the
constant $M$ estimates the interpolation norm of the family $\{T_w\}_w$.
Thus, in particular, for any $z\in K$, we have
$\|T_z(f_i(z))\|_z\le 2MD_K|\delta(z,w_0)|$.

At this stage, fix $z\in K$ and (using the hypothesis) pick
$\tilde{x}_i\in X_z$ with $\|\tilde{x}_i\|_z\le 2\kappa MC(K)|\delta(z,w_0)|$
and $T_z\tilde{x}_i=T_zf_i(z)$. Consequently, $\tilde{x}_i-f_i(z)$ belongs to
$\text{Ker}\,T_z$. If this latter space has dimension $<n$ then there exist
 scalars
$a_1,\ldots,a_n$ with $\max|a_i|=1$ and $\sum a_i\tilde{x}_i=\sum a_if_i(z)$.
Note that
$$
\left\|\sum a_if_i(z)\right\|_z\leq 2\kappa MD_KC(n,\rho({\Cal
F}))|\delta(z,w_0)|
$$
where $C(n,\rho({\Cal F}))$ is some constant depending exclusively on $n$ and
$\rho({\Cal F})$ (we may take in fact $C(n,\rho({\Cal F})):=n^{1/r}$
if ${\Cal F}$ is $r$-normed).
It then follows that there exists $g\in\Cal F$ with
$$
g(z)=\sum a_if_i(z)\,\,\,\text{and}\,\,\,
\|g\|_{\Cal F}\le 4\kappa MD_KC(n,\rho({\Cal F}))|\delta(z,w_0)|.
$$
Going further, the function $h:=(g-\sum a_if_i)\delta(\cdot,z)^{-1}$
belongs to ${\Cal F}$
and
$$
\|h\|_{\Cal F}\leq D_K\rho({\Cal F})\left(\|g\|_{\Cal F}
+\left\|\sum a_if_i\right\|_{\Cal F}\right)\leq C_1<+\infty
$$
where $C_1=C_1(n,K,M,\rho({\Cal F}))$ is some positive, finite constant.
In turn, this implies that
$$
\align
1=\text{max}\,|a_i|&\leq\left\|\sum a_ix_i\right\|_{w_0}=\left\|\sum
 a_if_i(w_0)\right\|_{w_0}\\
&\leq\rho({\Cal F})\left(\|g(w_0)\|_{w_0}+\left\|g(w_0)-\sum
 a_if_i(w_0)\right\|_{w_0}\right)\\
&\leq\rho({\Cal F})(\|g\|_{\Cal F}+|\delta(w_0,z)|\|h\|_{\Cal F})\\
&\leq C_2(|\delta(z,w_0)|+|\delta(w_0,z)|),
\endalign
$$
where $C_2=C_2(n,K,M,\rho({\Cal F}))$ is some positive, finite constant,
independent of $z$, $w_0\in K$. For $|\delta(z,w_0)|+|\delta(w_0,z)|$ small
enough this leads to a contradiction, thus the lemma is proved.
\hfill$\blacksquare$
\enddemo

The precise statement of the
Fredholm stability result alluded to earlier is contained in the next theorem.
Before we state it, we make one more definition. Consider $(U,\delta,Z,{\Cal
F})$
an interpolation scale and for each closed subspace $\tilde{Z}$ of $Z$ set
${\Cal F}_{\tilde{Z}}:=\{f\in{\Cal F};\,f(U)\subseteq\tilde{Z}\}$. Note that if
all admissible functions are continuous, then $(U,\delta,\tilde{Z},{\Cal
F}_{\tilde{Z}})$
becomes an interpolation scale in a natural way.

An interpolation scale as above is said to have the
{\it intersection property} provided there exists a
closed subspace $\tilde{Z}$ of $Z$ so that ${\Cal F}_{\tilde{Z}}\subseteq{\Cal
F}$
densely, and for each closed subspace
$E$ of $\tilde{Z}$ one has $E=\cap_{w\in U}\{f(w);\,f\in{\Cal F}_E\}$.

A couple of comments are in order here. The above definition formalizes the
idea
 of
an interpolation scale whose intersection of all intermediate spaces is fairly
 rich.
Indeed, it is easily checked that $\tilde{Z}$, referred to in the sequel as the
 {\it
intersection space}, has $\tilde{Z}\subseteq\cap_{w\in U}X_w$, and that in fact
$\tilde{Z}$ embeds into each $X_w$ continuously and densely.
Furthermore, the intersection property is not so difficult to check in many
 concrete
applications and, in fact,
in the case in which the interpolation scale arises by using the
classical complex method for a pair of compatible Banach spaces in the usual
 fashion,
then the intersection property is automatically satisfied.

\proclaim{Theorem 2.9} Let $\{T_w\}_{w\in U}$ be an interpolating  family
of operators between two interpolation scales $\{X_w\}_{w\in U}$,
$\{Y_w\}_{w\in U}$,
and suppose that $\{Y_w\}_{w\in U}$ has the intersection property.
Then the set of points $z\in U$ at which $T_z:X_z\to Y_z$ is a
Fredholm operator is open in $U$ and the index is locally constant on this set.

Moreover, for each $w\in U$ such that $T_w:X_w\to Y_w$ is Fredholm
there exists a neighborhood $\omega\subseteq U$ of $w$ such that
$\text{dim}\,\text{Ker}\,T_z\leq\text{dim}\,\text{Ker}\,T_w$
and $\text{dim}\,\text{Coker}\,T_z\leq\text{dim}\,\text{Coker}\,T_w$ for each
$z\in\omega$.
\endproclaim

\demo{Proof} To fix notation, recall we are assuming that
$\{Y_w\}_{w\in U}\equiv(U,\delta,W,{\Cal G})$ has the intersection property,
and
we let $\tilde{W}$ denote the intersection space. Also, let $Z$ and ${\Cal F}$
be,
respectively, the ambient space and the class of admissible functions for the
interpolation scale $\{X_w\}_{w\in U}$.

We claim that if $T_{w_0}$ is Fredholm for some $w_0\in U$ then $T_z$ is
Fredholm with the same index for $|\delta(z,w_0)|$ small enough.
To this effect, pick a finite dimensional subspace $E$ of $\tilde{W}$ so that
$E\oplus T_{w_0}(X_{w_0})=Y_{w_0}$. Augmenting $Z$ to $Z\oplus E$ and $\Cal F$
to
${\Cal F}\oplus{\Cal F}_E$ it is easy to check that
$(U,\delta,Z\oplus E,{\Cal F}\oplus{\Cal F}_E)$ is an interpolation scale
whose intermediate spaces are $\{X_w\oplus E\}_{w\in U}$.
Moreover, introducing $T'_w(x,y):=T_wx+y$, then $\{T'_w\}_{w\in U}$ becomes
an interpolating  family of operators between
the scales $\{X_w\oplus E\}_{w\in U}$ and $\{Y_w\}_{w\in U}$.

Since $T'_{w_0}:X_{w_0}\oplus E\to Y_{w_0}$ is onto, by Lemma 2.8,
the dimension of the kernel is constant in a neighborhood of $w_0$ where
$T'_{w_0}$ remains onto. This implies that $\dim T_z^{-1}(E)$ is constant
in this neighborhood. Now $\dim (\Coker T_z)= \dim E-\dim
T_z(X_z)\cap E=\dim E-\dim T_z^{-1}(E)+\dim (\Ker T_z)$ so that the index
is also constant on this neighborhood. Finally, note that also $\dim (\Coker
T_z)\le
\dim E=\dim(\Coker T_{w_0})$, and the conclusion follows.
\hfill$\blacksquare$
\enddemo

We conclude this section with a global stability result.

\proclaim{Theorem 2.10} Let $\{T_w\}_{w\in U}$ be an interpolating
family of operators between two interpolation scales, and suppose that
$U$, the space of parameters is connected and that $\{Y_w\}_{w\in U}$ has
the intersection property. Also, assume that there exists a point $w\in U$
such that $T_w:X_w\to Y_w$ is an isomorphism.

Then, if $\eta(T_z)>0$ for all $z\in U$, or if $\kappa(T_z)<\infty$ for all
$z\in U$ it follows that $T_z:X_z\to Y_z$ is an isomorphism for all $z\in U$.
\endproclaim

\demo{Proof} Let $V$ be set of $z$ where $T_z$ is an isomorphism.
By Theorem 2.7 $V$ is open. We show that $V$ is also closed.
Let $v\in V$ be a boundary point and assume that $\eta(T_v)>0$ but $T_v$ is not
onto.
Then we can as in the previous proof augment $Z$, the ambient space for the
interpolation scale $\{X_w\}_{w\in U}$ by some fixed, one-dimensional subspace
of the intersection space so that $T'_v$ still has a
lower bound $\eta(T'_v)>0$.
However, since $T'_z$ cannot be one-one for $z\in V$,
this will contradict Theorem 2.5.

For the other case assume $T_v$ onto but not one-one. Then a similar reasoning,
this time based on Lemma 2.8, leads to a contradiction and this proves the
result.
\hfill$\blacksquare$
\enddemo

\heading\S 3.\quad The complex method for pairs of quasi-Banach spaces
\endheading

The setup discussed in the previous section encompasses and unifies both
complex and real interpolation in a way which
avoids excessive discussion of technicalities.

Real interpolation, i.e. $(\theta,p)$-methods, can be obtained by starting from
 two
quasi-Banach spaces $X_0$, $X_1$, continuously embedded in some larger
 topological
vector space, and such that $X_0\cap X_1$ is dense in $X_j$, $j=0,1$.
We take $Z:=X_0+X_1$ to be the ambient space, the annulus $U:=\{z\in{\bold
 C};\,1<|z|<e\}$
as the space of parameters, and $\delta(z,w):=z-w$. Also, the class ${\Cal F}$
 of
admissible functions consists of analytic functions of the form
$f(z)=\sum_{n\in\bold Z}x_nz^n$, where $x_n\in X_0\cap X_1$, $n\in\bold Z$, and
the series converges in  the interior of the annulus. This latter space is
 equipped with
the quasi-norm given by
$\|f\|_{\Cal F}:=\text{max}\,\{(\sum_{n\in\bold Z}e^{n}\|x_n\|_{X_1}^p)^{1/p},
(\sum_{n\in\bold Z}\|x_n\|_{X_0}^p)^{1/p}\}$. Finally, the intermediate spaces
$\{X_\theta\}_{0<\theta<1}$ are as in (2.2) when $w=e^{\theta}$, $0<\theta<1$.
For more details we refer to the forthcoming paper \cite{12}.

Next we present the complex method of interpolation for pairs of quasi-Banach
 spaces.
Let us first recall that if $Y$ is a topological vector space
and
$U$ is an
open subset of the complex plane then a map $f:U\to Y$ is called analytic
if given $z_0\in U$ there exists $\eta>0$ so that there is a power
series expansion $f(z)=\sum_{n=0}^{\infty}y_nz^n$ converging uniformly
for $|z-z_0|<\eta$.
The theory of analytic functions with values in quasi-Banach spaces was
developed in \cite{33}, \cite{24}, \cite{25}.  It is important
to note here that in general there is no analogue of the Maximum Modulus
Principle, but there is an important subclass of quasi-Banach spaces
called A-convex in \cite{25} in which the Maximum Modulus
Principle does hold.

A quasi-Banach space $X$ is called {\it A-convex} if there is a constant
$C$ so that for every polynomial
$P:\bold C\to X$ we have $\|P(0)\|_X \le C\max_{|z|=1}\|P(z)\|_X$.  It
is shown in \cite{25} that if $X$ is A-convex it has an equivalent
quasi-norm which is plurisubharmonic (i.e. we can insist that $C=1$).
Let us also point out that being A-convex also implies that
$\max_{|z|\leq 1}\|f(z)\|_X\leq C\max_{|z|=1}\|f(z)\|_X$ for any analytic
function
$f:\{z;\,|z|<1\}\to X$ which is continuous on $\{z;\,|z|\leq 1\}$.

Next, we recall some basic properties of analytic functions:

\proclaim{Proposition 3.1}Suppose $0<p\le 1$ and that $m\in{\bold N}$ is
such that $m>\frac1p$.  Then there is a constant $C=C(m,p)$ so that if
$X$ is a $p$-normed quasi-Banach space and $f:\bar{\Delta}\to X$ is a
continuous function which is analytic on the unit disk $\Delta:=\{z;\,|z|<1\}$
then for $z\in \Delta$
we have $f(z)=\sum_{n=0}^{\infty}\frac{f^{(n)}(0)}{n!}z^n$,
and
$$
\|f^{(n)}(0)\|_X \le C(m+n)!\sup_{z\in \Delta}\|f(z)\|_X.
$$
\endproclaim

\noindent This is Theorem 6.1 of \cite{24}.

\proclaim{Proposition 3.2}Let $X$ be a quasi-Banach space and let $U$ be
an open subset of the complex plane.  Let $f_n:U\to X$ be a sequence of
analytic functions. If $\lim_{n\to\infty}f_n(z)=f(z)$ uniformly on
compacta then $f$ is also analytic.
\endproclaim

\noindent This follows from Theorem 6.3 of \cite{24}.

Our next theorem shows that, many times in practice, the ambient space $Z$
plays only a minor role in the setup. More specifically, assume that
$Z$ is a space of distributions in which a quasi-Banch space $X$ is
continuously
embedded. Then, having a $X$-valued function analytic for the quasinorm
topology
is basically the same as requiring analyticity for the weak topology
(induced on $X$ from $Z$).

\proclaim{Theorem 3.3} Suppose $X$ is a quasi-Banach space and that $U$ is
an open subset of the complex plane. Let $f:U\to X$ be a locally bounded
function. Suppose there is a weaker Hausdorff vector topology $\tau_0$
on $X$ which is locally $p$-convex for some $0<p<1$ and such that $f:U\to
(X,\tau_0)$ is analytic. Then $f:U\to X$ is analytic.
\endproclaim

\demo{Proof}We may suppose that $p$ is such that $X$ is also
$p$-normable.  Consider the largest locally $p$-convex vector topology
$\tau$  on $X$ which is weaker than the given quasi-norm topology and  so
that $f:U\to (X,\tau)$ is analytic. Next consider the
quasi-normed topology $\tau^*$ on $X$ generated by the unit ball which is
the $\tau-$closure of the given unit ball $B_X$ of $X$. It is easy to see
that $\tau^*$ is given by the quasi-norm $\|x\|_*=\sup_{j\in J}\|x\|_j$
where $\{\|\cdot\|_j:\ j\in J\}$ is the family of $\tau-$continuous
$p$-seminorms on $X$ satisfying $\|x\|_j\le \|x\|_X$.

We now show that $f:U\to (X,\tau^*)$ is analytic.  To this end suppose
$z_0\in U$.  Choose $r>0$ so that $\{z:|z-z_0|\le r\}\subset U$ and
$\sup_{|z-z_0|\le r}\|f(z)\|_X=M<\infty$. Then for each $i$ we can use
Proposition 3.1 to estimate $\|f^{(n)}(z_0)\|_j\le C(m+n)!r^{-n}M$ for
suitable $m>\frac1p$ and $C=C(m,p)$. But then $\|f^{(n)}(z_0)\|_*\le
C(m+n)!r^{-n}M^n$ and it follows that the Taylor series
$\sum_{n=0}^{\infty}\frac{f^{(n)}(z_0)}{n!}(z-z_0)^n$ converges in
$(X,\tau^*)$ necessarily to $f(z)$ for $|z-z_0|<r$. This proves our
assertion.

Now we must have $\tau=\tau^*$ by the choice of $\tau$.  This means
exactly that the identity map $i:X\to (X,\tau)$ is open (see
\cite{26} Theorem 1.4) by a form of the Open Mapping theorem. Hence
$\tau$ is the quasi-norm topology on $X$ and the result is
proved.
\hfill$\blacksquare$
\enddemo

We now let $(X_0,X_1)$ be an interpolation couple of
quasi-Banach spaces, i.e. $X_j$, $j=0,1$ are
continuously embedded into a larger topological vector space $Y$, and
$X_0\cap X_1$ is dense in $X_j$, $j=0,1$. Also, let $U$
stand in this case for the strip
$\{z\in{\bold C};\,0<\text{Re}\,z<1\}$, and take $Z:=X_0+X_1$,
$\delta(z,w):=z-w$.

We then define ${\Cal F}$, the class of admissible functions as the space
of bounded, analytic functions $f:U\to X_0+X_1$ which extend
continuously to the closure of the strip such that the traces
$t\mapsto f(j+it)$ are bounded continuous functions into $X_j$, $j=0,1$.

As in \cite{22}, we endow ${\Cal F}$ with the quasi-norm
$$
\align
\|f\|_{\Cal
 F}:=&\text{max}\left\{\underset{t}\to{\text{sup}}\,\|f(it)\|_{X_0}\,,\,
\underset{t}\to{\text{sup}}\,\|f(1+it)\|_{X_1},
\sup_{z\in U}\|f(z)\|_{X_0+X_1}\right\}.
\tag{3.1}
\endalign
$$
With this definition, the axioms (1)-(3) stated in section 2
are readily verified, and we denote
by $X_w=[X_0,X_1]_w$, $w\in U$, the corresponding intermediate spaces.
Of course $X_w$ depends only on $\Re w$ and so we typically consider
$X_{\theta}$ for $0<\theta<1$.
Note that the second term in the definition of $\|\cdot\|_{\Cal F}$ is
precisely
designed to ensure both the completeness of ${\Cal F}$ and the continuity of
the evaluation functions. However, it should be pointed out that this term can
be dropped as being dominated by the first if $X_0+X_1$ is A-convex.

Note that if $T$ is a linear bounded operator mapping $X_j$
into itself,
$j=0,1$. Then $T$ is clearly bounded on $X_0+X_1$ and is therefore an
interpolating operator. Interpolating analytic families of operators
$\{T_w\}_{0<\Re w<1}$ can also be defined as in the previous section.

Our definition which is the natural extension of Calder\'on's
original
definition for Banach spaces is not, however, usually employed in the
existing literature.  In order to avoid discussing analytic functions in
quasi-Banach spaces, authors have typically employed two variants.

One variant of this approach is to replace $\Cal F$ by its closed
subspace $\Cal F_0$ generated by functions with only
finite-dimensional range, and this is used in many places in
the literature e.g. \cite{22},\cite{21} and \cite{13}.
While it is not clear in general that leads to the same interpolation
spaces, in the case we wish to consider (Hardy spaces) this makes no
essential difference.  The advantage here is simply that one can use
functions with finite-dimensional range, without discussing the meaning
of analyticity.

A second variant is to require that the functions $f\in\Cal F$ be
analytic into the ambient space $Z$, which is typically a locally convex
space of distributions (cf. \cite{22}, \cite{6}). Theorem 3.3 shows
that this does not change $\Cal F$; however if one requires only
$Z$-continuity at the boundary this may result in a different space.
In \cite{2} Bernal takes $Z$ to be A-convex for certain applications.

We next point out that the method described above gives the result
predicted by the Calder\'on formula for nice pairs of function
spaces. The theorem below is due to Gomez and Milman \cite{21} for a
$\Cal F_0-$variation of the definition. We will give the details
for the convenience of the reader.

Let $(\Omega,\mu)$ be a $\sigma-$finite measure space and let
$L_0$ be the collection of all complex-valued, $\mu$-measurable
functions on $\Omega$. Recall
that {\it a quasi-Banach function space} $X$ on $(\Omega,\mu)$ is
an order-ideal in the space $L_0$ containing
a strictly positive function, equipped with a quasi-norm $\|\cdot\|_X$
so that $(X,\|\cdot\|_X)$ is complete, and if
$f\in X$ and $g\in L_0$ with $|g|\le |f|$ a.e. then $g\in X$ with
$\|g\|_X\le \|f\|_X$.

\proclaim {Theorem 3.4}Let $\Omega$ be a Polish space and
let $\mu$  be a
$\sigma-$finite  Borel
measure on $\Omega$.  Let $X_0,X_1$ be a pair of quasi-Banach function
spaces on $(\Omega,\mu)$.  Suppose that both $X_0$ and $X_1$ are
A-convex and separable.
Then $X_0+X_1$ is A-convex and
$X_{\theta}=[X_0,X_1]_{\theta}=X_0^{1-\theta}X_1^{\theta},$ in the sense
of equivalence of quasi-norms.
\endproclaim

\demo{Remarks}$X_0^{1-\theta}X_1^{\theta}$ is the function space
defined by the quasi-norm
$$
\|h\|:=\inf\{\|f\|_{X_0}^{1-\theta}\|g\|_{X_1}^{\theta};\ |h|\le
|f|^{1-\theta}|g|^{\theta},\,\,f\in X_0,\,\,g\in X_1\}.
$$
The hypothesis of separability in this case is equivalent to
order-continuity, which can reformulated as the following property: if
$g\in X_j,\
(j=0,1)$ and
$|f_n|\le
|g|$ for all $n$ and $f_n\to f$ a.e. then $\|f_n-f\|_{X_j}\to 0$.
\hfill$\square$
\enddemo

\demo{Proof}We indicate briefly the argument that $X_0+X_1$ is A-convex.
We may assume (\cite{25}) that $X_0$ and $X_1$ are both $p$-convex
lattices i.e. we have for $f_1,\ldots,f_n\in X_j$ we have
$$
\left\|(\sum_{k=1}^n|f_k|^p)^{1/p}\right\|_{X_j}\le
\left(\sum_{k=1}^n\|f_k\|_{X_j}^p\right)^{1/p}
$$
for $j=0,1$. This implies that the spaces
$X_j^p:=\{f\in L_0; |f|^{1/p}\in X_j\}$ when normed by
$\|f\|_{X_j^p}=\||f|^{1/p}\|_{X_j}^p$ are Banach function spaces.
Now suppose $f_k\in X_0+X_1=W$ for $k=1,2,\ldots,n$.
Pick $g_k\in X_0$ and $h_k\in X_1$ with $\|g_k\|_{X_0}+\|h_k\|_{X_1}\le
2\|f_k\|_{W}$. Then
$$
\align
\left\|\left(\sum_{k=1}^n|f_k|^p\right)^{1/p}\right\|_W &\le
\left\|\left(\sum_{k=1}^n(|g_k|+|h_k|)^p\right)^{1/p}\right\|_W\\
&\le
2^{1/p-1}\left\|\left(\sum_{k=1}^n|g_k|^p\right)^{1/p}+
\left(\sum_{k=1}^n|h_k|^p\right)^{1/p}\right\|_W \\
&\le 2^{1/p-1}\left\|\left(\sum_{k=1}^n|g_k|^p\right)^{1/p}\right\|_{X_0} +
2^{1/p-1}\left\|\left(\sum_{k=1}^n|h_k|^p\right)^{1/p}\right\|_{X_1}\\
&\le
2^{1/p-1}\left(\sum_{k=1}^n\|g_k\|_{X_0}^p\right)^{1/p}+2^{1/p-1}
\left(\sum_{k=1}^n\|h_k\|_{X_1}^p\right)^{1/p}\\
&\le  2^{1/p+1}\left(\sum_{k=1}^n\|f_k\|_W^p\right)^{1/p}
\endalign
$$
so that $X_0+X_1$ is also a $p$-convex lattice and so is A-convex.
As in \cite{21} we remark that the remainder of the argument is very
similar to the original argument of A. P. Calder\'on (\cite{5}).

First suppose $f_0\in X_0$ and $f_1\in X_1$ and suppose that $f_0$
and $f_1$ have the same supports.  Then the function $z\to
|f_0|^{1-z}|f_1|^z$ is easily shown to be in $\Cal F$ (here we use the
separability hypotheses on $X_0$ and $X_1$, and it follows that
$$
\||f_0|^{1-\theta}|f_1|^{\theta}\|_{\theta}\le
C\| f_0\|_{X_0}^{1-\theta}\|f_1\|_{X_1}^{\theta}.
$$
for a suitable constant $C$.  This implies that for any $f\in
[X_0,X_1]_{\theta}$ we have $\|f\|_{\theta}\le
C\|f\|_{X_0^{1-\theta}X_1^{1-\theta}}$.
(Here $C=1$ if we drop the third term in our definition of $\Cal F$.)

Now suppose $F\in \Cal F$. For any $q\le p$ let $W_q=(X_0+X_1)^{q}$
which is also a Banach function space. We first observe that if $z_0\in
U$ then for some $r>0$ we have a series expansion
$F(z)=\sum_{n=0}^{\infty}f_n(z-z_0)^n$ with $f_n\in  X_0+X_1$ valid
for $|z-z_0|<r$ and such that
$\limsup\|f_n\|_{X_0+X_1}^{1/n}\le r^{-1}$. Then $\sum_{n=0}^{\infty}
|f_n||z-z_0|^n$ is convergent in $X_0+X_1$ and it follows easily that we
can argue pointwise to deduce that
$$
|F(z_0)(\omega)|^q\le\frac{1}{2\pi}
\int_0^{2\pi}|F(z_0+\rho e^{it})(\omega)|^q\,dt
$$
for $\mu$-a.e. $\omega\in\Omega$ and any $\rho<r$. Now $z\mapsto |F(z)|^q$ is
continuous into $W_q$ and a simple application of Fubini's theorem shows
that for any positive $\phi\in W_q^*$ we have $z\mapsto\phi(|F(z)|^q)$ is
subharmonic on $U$. Thus
$$
\phi(|F(\theta)|^q)\le\int_{-\infty}^{\infty}P_0(\theta,t)\phi(|F(it)|^q)dt
+\int_{-\infty}^{\infty}P_1(\theta,t)\phi(|F(1+it)|^q)dt
$$
where $P_0$ and $P_1$ are the components of the Poisson kernel for the
strip. Here $\int P_0dt =1-\theta$ and $\int P_1\,dt=\theta$.
Let
$g=((1-\theta)^{-1}\int_{-\infty}^{\infty}P_0(\theta,t)|F(it)|^pdt)^{1/p}$.
By
$p$-convexity (or by the local convexity of $W_p$) we have $g\in X_0$ and
$\|g\|_{X_0}\le \|F\|_{\Cal F}$; similarly if
$h=(\theta^{-1}\int_{-\infty}^{\infty}P_1(\theta,t)|F(1+it)|^pdt)^{1/p}$
then
$\|h\|_{X_1}\le \|F\|_{\Cal F}$. Now if $q\le p$ we deduce $|F(\theta)|^q
\le (1-\theta)g^q+\theta h^q$. Taking $log$ of both sides and then
letting $q\to 0$, this implies $|F(\theta)|\le g^{1-\theta}h^{\theta}$ or
$\|F(\theta)\|_{X_0^{1-\theta}X_1^{\theta}}\le \|F\|_{\Cal F}$.
\hfill$\blacksquare$
\enddemo

For $0<p<\infty$ let $H^p=H^p({\bold R}^m)$ denote the real variables Hardy
spaces as introduced in \cite{18}.
That the family of Hardy spaces $\{H^p({\bold R}^m)\}_{0<p<\infty}$ is an
interpolation scale of A-convex spaces for the complex method is essentially
well known (cf., e.g., \cite{6}, \cite{13}, \cite{19}, \cite{1}), although
the precise definition of complex interpolation can vary.

For our purposes, it is useful to note that a simple way to see
that
$[H^p,H^q]_{\theta}=H^r$ where
$\frac1r=\frac{1-\theta}{p}+\frac{\theta}{q}$ would be to use the
identification with a pair of lattices given by taking a common
unconditional basis in the spaces $H_p,H_q$.  This can be achieved by
using a wavelet basis (cf. Theorem 7.20 of \cite{20}) and reduces the
problem to interpolating two sequence spaces $\dot{f}^{0,2}_p$ and
$\dot{f}^{0,2}_q$. Now Theorem
3.4
(with
$\Omega=\bold N$ and $\mu$ counting measure) can be applied; we need
only the appropriate Calder\'on formula for the $f$-spaces and this is
and done in \cite{19}. A rather similar approach is given via tent
spaces in \cite{1}.

\heading\S 4.\quad Applications to PDE's
\endheading

Let $\Omega$ be a Lipschitz domain in ${\bold R}^m$. In the unbounded case this
means that $\Omega$ is the domain above the graph of a Lipschitz function
$\varphi:{\bold R}^{m-1}\to{\bold R}$, i.e.
$\Omega=\{(\varphi(x)+t,x);\,x\in{\bold R}^{m-1},\,\,t>0\}$. Also, in the
 bounded case,
$\Omega$ is a bounded domain whose boundary is locally given by graphs of
 Lipschitz
functions. Let $n(P)$ be the outward unit normal defined at almost every
boundary point $P\in\partial\Omega$ and denote by $d\sigma$ the surface measure
on $\partial\Omega$. Set ${\sigma}_m$ for the area of the unit sphere
$S^{m-1}\subseteq{\bold R}^m$ and $\langle\cdot,\cdot\rangle$ for the usual
inner product of vectors in ${\bold R}^m$ (the corresponding Euclidean norm in
${\bold R}^m$ is denoted by $|\cdot|$). We denote by $L^p(\partial\Omega)$ the
Lebesgue spaces of scalar-valued functions which are
$d\sigma$-measurable and $p$-th power integrable on $\partial\Omega$.

Next, recall that the classical {\it double layer potential operator} on
 $\partial\Omega$
acts on boundary densities $f:\partial\Omega\to\bold C$ by
$$
{\Cal D}
 f(X):=\frac{1}{\sigma_m}\int_{\partial\Omega}\frac{\langle{n(Q)},{Q-X}\rangle}
{|X-Q|^m}f(Q)\,d\sigma(Q),\quad X\in\Omega.
\tag{4.1}
$$
Then, at almost every boundary point $P\in\partial\Omega$ we have the usual
jump relations
$$
\align
\lim_{\Sb X\to P\\X\in\gamma(P)\endSb}{\Cal D}f(X)&=\frac{1}{2}f(P)+
\text{p.v.}\frac{1}{\sigma_m}
\int_{\partial\Omega}\frac{\langle{n(Q)},{Q-P}\rangle
 }
{|P-Q|^m}f(Q)\,d\sigma(Q)\\
&=:\left(\frac{1}{2}I+{\Cal K}\right)f(P).
\tag{4.2}
\endalign
$$
Here, $\gamma(P):=
\{X\in\Omega ;\,|X-P|\leq 2\,\text{dist}\,(X,\partial\Omega)\}$
is the nontangential approach region corresponding to $P\in\partial\Omega$.
Also, $I$ is the identity operator and the principal value
singular integral ${\Cal K}f$ is the so-called {\it singular double layer
 potential
operator} acting on $f$.

We shall also be interested to work with ${\Cal K}^*$, the formal
adjoint of the singular double layer potential operator, defined by
$$
{\Cal K}^*f(P)=
\text{p.v.}\frac{1}{\sigma_m}
\int_{\partial\Omega}\frac{\langle{n(P)},{P-Q}\rangle
 }
{|P-Q|^m}f(Q)\,d\sigma(Q),\qquad P\in\partial\Omega .
$$
More material on these can be found in e.g. \cite{17}, \cite{23},
\cite{35}, \cite{14}.

For each $0<\alpha<1$, denote by $C^{\alpha}(\partial\Omega)$ the class of
 H\"older
continuous functions of order $\alpha$ on $\partial\Omega$, i.e.
$f\in C^{\alpha}(\partial\Omega)$ if $|f(P)-f(Q)|\leq M\,|P-Q|^{\alpha}$ for
some finite constant $M$, uniformly for $P$, $Q\in\partial\Omega$. Set
$\|f\|_{C^{\alpha}(\partial\Omega)}$ for the best constant in the previous
 inequality.

Our first result in this section is the integral representation formula
 contained
in the following theorem.

\proclaim{Theorem 4.1} There exits a small, positive $\alpha_0$ which depends
 only on the
dimension $m$ and the Lipschitz character of the domain $\Omega$ for which the
following property holds. If a function $u$ harmonic in $\Omega$ (in the
 unbounded
domain case $u$ is also assumed to be suitably small at infinity)
assumes continuously boundary values and
$u|_{\partial\Omega}\in C^{\alpha}(\partial\Omega)$ with $0<\alpha<\alpha_0$,
 then $u$
is the double layer extension in $\Omega$ of a boundary density from the same
 class.
In other words, $u$ can be represented in the form
$$
u(X)=\frac{1}{\sigma_m}\int_{\partial\Omega}\frac{\langle{n(Q)},{Q-X}\rangle}
{|X-Q|^m}f(Q)\,d\sigma(Q),\quad X\in\Omega,
$$
for some boundary density $f$ which is
H\"older continuous of order $\alpha$ on $\partial\Omega$.
Moreover, $f$ is uniquely determined and in fact
$f=(\frac{1}{2}I+{\Cal K})^{-1}(u|_{\partial\Omega})$.
\endproclaim

Note that an immediate corollary is the classical estimate
$\|u\|_{C^\alpha(\bar{\Omega})}\leq C\|u\|_{C^\alpha(\partial\Omega)}$
valid for any function $u$ harmonic in $\Omega$, provided $\alpha>0$ is
sufficiently small, depending on $\partial\Omega$.
Before proceeding any further, let us also point out that similar results
are valid for solutions to the Dirichlet problem with H\"older continuous
boundary data for the Helmholtz operator $\triangle +k^2$ in Lipschitz
domains. Of course, here the wave number $k\in{\bold C}$, $\text{Im}\,k\geq 0$,
is assumed not to be a Dirichlet eigenvalue for the domain under discussion.

Recall that for $\frac{m-1}{m}<p\leq 1$, $p$-atoms are defined as
scalar-valued functions $a$ supported in a surface
ball $S(P,r):=B(P,r)\cap\partial\Omega$
(for some $P\in\partial\Omega$ and some $r>0$), such that
$\|a\|_{L^\infty}\leq\sigma(S)^{-\frac{1}{p}}$ and
 $\int_{\partial\Omega}a\,d\sigma =0$.
Set $\alpha:=\frac{1-p}{p}(m-1)\in(0,1)$ so that the pairing
$\int_{\partial\Omega}a\,h\,d\sigma$ between a $p$-atom $a$ and a function $h$
H\"older continuous of order $\alpha$ on $\partial\Omega$ is well defined.
Then the atomic Hardy space $H^p_{\text{at}}(\partial\Omega)$ (cf.
\cite{10}) is defined as the collection of all continuous functionals
$f$ on $C^{\alpha}(\partial\Omega)$  which can be represented in the form
$f=\sum_j\lambda_ja_j$, where $a_j$'s are $p$-atoms and
$\lambda_j\in{\bold R}$ are so that $\sum_j|\lambda_j|^p<+\infty$.
It is well known that this space is complete when endowed with the natural
 quasi-norm
$$
\|f\|_{H^p_{\text{at}}(\partial\Omega)}
:=\text{inf}\,\left\{\left(\sum_j|\lambda_j|^p\right)^{\frac{1}{p}};
\,f= \sum_j\lambda_ja_j,\,\,\lambda_j\in{\bold
 R},\,\,a_j\,\,p-\text{atom}\right\}
$$
(in fact, the above formula defines a norm if and only if $p=1$).
Accordingly, the well-definiteness of the pairing mention above translates into
$$
(H^p_{\text{at}}(\partial\Omega))^*=C^{\alpha}(\partial\Omega).
\tag{4.3}
$$

Next, we introduce the scale $\{H^p(\partial\Omega)\}_{\frac{m-1}{m}<p<\infty}$
 by
$$
H^p(\partial\Omega):=
\cases
L^p(\partial\Omega),\,\,\text{if}\,\,1<p<+\infty,\\
H^p_{\text{at}}(\partial\Omega),\,\,\text{if}\,\,\frac{m-1}{m}<p\leq 1.
\endcases
$$

A simple observation which is of importance for us is that the $\{H^p\}_p$
scale is essentially invariant under bi-Lipschitz changes of coordinates.
More concretely, if $D$ is another Lipschitz domain in ${\bold R}^m$
and $\Lambda:\partial D\to\partial\Omega$ is a bi-Lipschitz homeomorphism then,
for each $\frac{m-1}{m}<p<+\infty$,
$$
H^p_{\text{at}}(\partial\Omega)\ni f\mapsto (f\circ\Lambda)|J_{\Lambda}|\in
H^p_{\text{at}}(\partial D)
$$
is an isomorphism of quasi-Banach spaces. Here $J_{\Lambda}$ is the Jacobian of
$\Lambda$ (note that $|J_{\Lambda}|$ is bounded away from zero and infinity
 uniformly
on $\partial D$) and the idea is that $p$-atoms are mapped into a certain fixed
 multiple
of $p$-atoms.
This observation can be used to transfer problems formulated {\it entirely} on
 the scale
$\{H^p(\partial\Omega)\}_p$ to more convenient ones, like $\{H^p({\bold
 R}^{m-1})\}_p$
(in the case of unbounded Lipschitz domains) and $\{H^p(S^{m-1})\}_p$
(in the case of bounded Lipschitz domains).

The proof of Theorem 4.1 will be accomplished in a series of lemmas which we
 state below.
Fix some Lipschitz domain $\Omega$ in ${\bold R}^m$.
First, recall the double layer potential operator ${\Cal D}$ from (4.1) and
denote by $C^{\alpha}(\overline{\Omega})$ the space of all functions in
$\Omega$ which extend as H\"older continuous functions
of order $\alpha$ in $\overline{\Omega}$ endowed with the usual homogeneous
norm
$\|u\|_{C^{\alpha}(\overline{\Omega})}
:=\text{sup}\,\{\frac{|u(X)-u(Y)|}{|X-Y|^{\alpha}};\,X\neq
 Y,\,\,X,Y\in\Omega\}$.

\proclaim{Lemma 4.2} For each $0<\alpha<1$,
${\Cal D}:C^{\alpha}(\partial\Omega)\to C^{\alpha}(\overline{\Omega})$
is well defined, linear and bounded.
\endproclaim

This is perhaps well known but, since we lack an exact reference, for the
sake of completeness we indicate a simple proof.

\demo{Proof} It is not too difficult to check that a continuously
differentiable
function $F:\Omega\to{\bold R}$ extends in $C^{\alpha}(\overline{\Omega})$
 provided
$$
\underset{X\in\Omega}
\to{\text{sup}}\,{\text{dist}\,(X,\partial\Omega)}^{1-\alpha
 }
|\triangledown F(X)|\leq C_0<+\infty ,
\tag{4.4}
$$
and that in fact
$\|F\|_{C^{\alpha}(\overline{\Omega})}\leq C_1C_0$, where $C_1$ depends only on
$\partial\Omega$, $m$ and $\alpha$. Next, pick $f\in
 C^{\alpha}(\partial\Omega)$,
$X\in\Omega$ and set $F:={\Cal D}f$, $h:=|X-P|$, where
$P\in\partial\Omega$ is such that $|X-P|=\text{dist}\,(X,\partial\Omega)$.
Note that, since ${\Cal D}1=1$, there is no loss of generality in assuming that
$f(P)=0$. Split the domain of integration in ${\Cal D}f$ over
$\{Q\in\partial\Omega;\,|P-Q|\leq 100h\}$ and $\{Q\in\partial\Omega;\,|P-Q|\geq
 100h\}$.
If we now bound the kernel of $\triangledown{\Cal D}$ by $Ch^m$
in the first resulting integral and by $C|P-Q|^m$ in the second one,
then straightforward size estimates yield (4.4).
\hfill$\blacksquare$
\enddemo

\proclaim{Lemma 4.3} For any $\frac{m-1}{m}<p<+\infty$, the operator ${\Cal
 K}^*$ is
a well defined, bounded linear mapping of $H^p(\partial\Omega)$ into itself.
\endproclaim

\demo{Proof} This is essentially well known too. The case $1<p<+\infty$ is
 basically
the theorem of Coifman, McIntosh and Meyer (eventually combined with the
method of rotation of Calder\'on and Zygmund); see e.g. \cite{9},
\cite{34}
for details. In turn, this can be used to show that ${\Cal K}^*$ maps
atoms into molecules so that case $\frac{m-1}{m}<p\leq 1$ is covered by the
rather general theory in \cite{10}.
\hfill$\blacksquare$
\enddemo

\proclaim{Lemma 4.4} The operator $\frac{1}{2}I+{\Cal K}^*$ is
an isomorphism of $H^1(\partial\Omega)$ onto itself.
\endproclaim

\demo{Proof} This is Theorem 3.9 in \cite{14}.
\hfill$\blacksquare$
\enddemo

We are now in a position to present the proof of Theorem 4.1.

\demo{Proof of Theorem 4.1} From the discussion in the last part of
section 3
we know that $\{H^p\}_{\frac{m-1}{m}<p<\infty}$ is an interpolation
scale and that $\frac{1}{2}I+{\Cal K}^*$ is an interpolation operator on it.
At this point, Lemma 4.3, Lemma 4.4, and
the discussion in section 3 together with the results in section 2
give us that there exists a small, positive $\epsilon$
such that $\frac{1}{2}I+{\Cal K}^*$ is an isomorphism of $H^p$ onto itself for
 each
$1-\epsilon<p\leq 1$. Now, by duality, there exists some small positive
 $\alpha_0$
such that $\frac{1}{2}I+{\Cal K}$ is an isomorphism of $C^{\alpha}$ onto itself
 for each
$0<\alpha<\alpha_0$ (as pointed out in the introduction, this has been first
observed in
\cite{4} with a completely different proof).
With this at hand, the conclusion in Theorem 4.1 easily  follows Lemma 4.2 and
the
jump relation (4.2).
\hfill$\blacksquare$
\enddemo

\demo{Remark} It is well known that
$\{C^{\alpha}(\partial\Omega)\}_{0<\alpha<1}$
is a complex interpolation scale. It then follows from the global stability
result
in section 2 (i.e. Theorem 2.10) that, with regard to the invertibility of
$\frac{1}{2}I+{\Cal K}$ on $C^{\alpha}(\partial\Omega)$, what breaks down for
$\alpha$
close to $1$ (on general Lipschitz domains) is precisely
the closedness of the range of this operator.

The same considerations apply, for instance, to the operator
$\frac{1}{2}I-{\Cal K}^*$
when acting on $L^p_0(\partial\Omega)$ for $p$ large ($p>2$).
\hfill$\square$
\enddemo

Recall next the {\it single layer potential operator} ${\Cal S}$ which
acts on boundary densities by
$$
{\Cal S}f(P):=\int_{\partial\Omega}\frac{1}{|P-Q|^{m-1}}f(Q)\,d\sigma(Q),
$$
and
let $\triangledown_{\text{tan}}:=\triangledown -n\frac{\partial}{\partial n}$
be the usual {\it tangential gradient operator} on $\partial\Omega$.
Following \cite{14},
we say that $a:\partial\Omega\to{\bold R}$ is a {\it regular atom} if
$\text{supp}\,a\subseteq S(P,r)$ for some $P\in\partial\Omega$, $r>0$ and
$\|\triangledown_{\text{tan}}a\|_{L^{\infty}(\partial\Omega)}\leq
\sigma(S)^{-\frac{1}{p}}$. Then define $H^{p,1}_{\text{at}}(\partial\Omega)$
as the $l^p$-span of regular atoms.

\proclaim{Theorem 4.5} Let $\Omega$ be a Lipschitz domain in ${\bold R}^m$.
Then
there exists a small, positive $\epsilon=\epsilon(\partial\Omega ,m)$ such that
the operator
${\Cal S}:H^p_{\text{at}}(\partial\Omega)\to
H^{p,1}_{\text{at}}(\partial\Omega)/{\bold R}$
is an isomorphism for each $1-\epsilon<p\leq 1$.
\endproclaim

\demo{Proof} The case $p=1$ is contained in \cite{14} and, hence, the
conclusion follows by the same pattern by proving that
$H^{p,1}_{\text{at}}(\partial\Omega)/{\bold R}$ is an interpolation scale for
the complex method. In fact, this last point can be circumvented as follows.
First, using the results in \cite{14} and Theorem 2.5 we see that there exists
$\epsilon >0$ such that
$$
\|f\|_{H^p_{\text{at}}(\partial\Omega)}\leq
C\|\triangledown_{\text{tan}}Sf\|_{H^p_{\text{at}}(\partial\Omega,{\bold R}^3)}
\leq
C\|Sf\|_{H^{p,1}_{\text{at}}(\partial\Omega,{\bold R}^3)}
$$
for any $f\in H^p_{\text{at}}(\partial\Omega)$, $1-\epsilon <p\leq 1$.
In particular, ${\Cal S}:H^p_{\text{at}}(\partial\Omega)\to
H^{p,1}_{\text{at}}(\partial\Omega)/{\bold R}$ is injective with closed range.
Now the conclusion follows from the fact that ${\Cal
S}:L^2_0(\partial\Omega)\to
L^{2,1}(\partial\Omega)/{\bold R}$ is an isomorphism; see \cite{34}.
\hfill$\blacksquare$
\enddemo

Our last application concerns the Lam\'e system of linear elastostatics on a
Lipschitz domain $\Omega$ in ${\bold R}^3$
$$
\mu\triangle\vec{u}+(\lambda
 +\mu)\triangledown(\text{div}\,\vec{u})=\vec{0}\,\,\,
\text{in}\,\,\,\Omega
\tag{4.5}
$$
where $\mu>0$ and $\lambda>-\frac{2}{3}\mu$ are the so-called Lam\'e constants.
Natural boundary conditions for (4.5) can be obtained by prescribing
$$
\vec{u}|_{\partial\Omega}=\vec{g},
\tag{4.6}
$$
(i.e. Dirichlet boundary conditions), or
$$
\left.\frac{\partial\vec{u}}{\partial\nu}:=(\lambda(\text{div}\,\vec{u})n
+\mu[\triangledown\vec{u}+\triangledown\vec{u}^t]n)\right|_{\partial\Omega}=\vec
 {h},
\tag{4.7}
$$
where the boundary conormal derivative operator $\partial/\partial\nu$ is
called
{\it traction} or {\it stress}. Note that the superscript $t$ indicates
 transpose of
the $3\times 3$ matrix $\triangledown\vec{u}=(\partial_iu_j)_{i,j}$.

Recall the Kelvin matrix $\Gamma=(\Gamma_{i,j})_{i,j}$ of fundamental solutions
 for
the system of elastostatics, where for $1\leq i,j\leq 3$
$$
\Gamma_{i,j}(X):=\frac{1}{2\omega_2}\left(\frac{1}{\mu}+\frac{1}{2\mu
 +\lambda}\right)
\frac{\delta_{i,j}}{|X|}+
\frac{1}{2\omega_2}\left(\frac{1}{\mu}-\frac{1}{2\mu +\lambda}\right)
\frac{X_iX_j}{|X|^3}.
$$
See \cite{27}. For vector valued densities on $\partial\Omega$ we define
the {\it single} and the {\it double} layer potential operators by
$$
{\bold
 S}\vec{f}(X):=\int_{\partial\Omega}\Gamma(X-Q)\,\vec{f}(Q)\,d\sigma(Q),\quad
 X\in
{\bold R}^3
$$
and by
$$
{\bold D}\vec{f}(X):=\int_{\partial\Omega}\left(\frac{\partial}{\partial\nu(Q)}
\Gamma(X-Q)\right)^t\,\vec{f}(Q)\,d\sigma(Q),\quad X\in\Omega ,
$$
respectively. Here the operator $\partial/\partial\nu$ applies to each column
of
the matrix $\Gamma$.

Finally, we record the corresponding jump relation
$$
\align
\lim_{\Sb X\to P\\X\in\gamma(P)\endSb}{\bold
 D}\vec{f}(X)&=\frac{1}{2}\vec{f}(P)+
\text{p.v.}\int_{\partial\Omega}\left(\frac{\partial}{\partial\nu(Q)}
\Gamma(X-Q)\right)^t\,\vec{f}(Q)\,d\sigma(Q)\\
&=:\left(\frac{1}{2}{\bold I}+{\bold K}\right)\vec{f}(P),\quad
 P\in\partial\Omega ,
\tag{4.8}
\endalign
$$
and let ${\bold K}^*$ denote the formal adjoint of ${\bold K}$. See \cite{16}
 for more
details.

\proclaim{Theorem 4.6} Let $\Omega$ be a Lipschitz domain in ${\bold R}^3$.
Then there exists a small, positive $\epsilon$
such that $\frac{1}{2}{\bold I}+{\bold K}^*$ is an isomorphism of
$H^p_{\text{at}}(\partial\Omega ,{\bold R}^3)$ onto itself for each
$1-\epsilon<p\leq 1$.

Also, there exists some small positive $\alpha_0$
such that $\frac{1}{2}{\bold I}+{\bold K}$ is an isomorphism of
$C^{\alpha}(\partial\Omega,{\bold R}^3)$ onto itself for each
 $0<\alpha<\alpha_0$.
\endproclaim

\demo{Proof} The fact that ${\bold K}^*$ maps
$H^p_{\text{at}}(\partial\Omega ,{\bold R}^3)$ boundedly into itself for each
$\frac{m-1}{m}<p\leq 1$ follows essentially as before (cf. Lemma 4.3). Also,
by combining Lemma 1.6 in \cite{15} with the techniques in \cite{14} one can
show that $\frac{1}{2}{\bold I}+{\bold K}^*$ is an isomorphism of
$H^1_{\text{at}}(\partial\Omega ,{\bold R}^3)$. Now the first part of the
 theorem
follows from this, the discussion in section 3 and the stability results in
 section 2.
The second part is a direct corollary of this and (4.3).
\hfill$\blacksquare$
\enddemo

The next theorem provides an alternative approach to results of B. Dahlberg and
C. Kenig (see \cite{15} for a proof based on different techniques). Let us
 point out
that the layer potential integral representation formula (4.9) is new.

\proclaim{Theorem 4.7} Let $\Omega$ be a Lipschitz domain in ${\bold R}^3$.
Then there exists some small positive $\alpha_0$ such that for each
 $0<\alpha<\alpha_0$
the boundary value problem (4.5)-(4.6) has a unique solution
$\vec{u}\in C^{\alpha}(\overline{\Omega},{\bold R}^3)$ for any boundary
data $\vec{g}\in C^{\alpha}(\partial\Omega,{\bold R}^3)$.

Moreover, we have the integral representation formula
$$
\vec{u}={\bold D}\vec{f}\quad\text{in}\,\,\Omega
\tag{4.9}
$$
with $\vec{f}:=(\frac{1}{2}{\bold I}+{\bold K})^{-1}\vec{g}\in
C^{\alpha}(\partial\Omega,{\bold R}^3)$, and the natural estimate
$\|\vec{u}\|_{C^{\alpha}(\overline{\Omega},{\bold R}^3)}\leq
C\|\vec{g}\|_{C^{\alpha}(\partial\Omega,{\bold R}^3)}$
holds for some $C=C(\alpha ,\partial\Omega)>0$.
\endproclaim

\demo{Proof} This is a direct corollary of Theorem 4.6 and the jump relation
(4.8) much in the spirit of the proof of Theorem 4.1
(note that the obvious analogue of Lemma 4.2 remains valid for the case we are
considering here).
\hfill$\blacksquare$
\enddemo

Further consequences of the above theorems in the study of elliptic
boundary value problems in non-smooth domains are analyzed in \cite{29}.
We conclude with a number of observations of independent interest.

\vskip 0.08in
\noindent{\bf Remarks.}

(1) Approximating continuous functions by H\"older continuous on the boundary
and relying on $L^{\infty}$ estimate from \cite{15}, it follows from Theorem
4.7 that
solutions $\vec{u}$ of (4.5)-(4.6) with $\vec{g}$ continuous are in fact
continuous up to and including the boundary of the domain.

(2) Let us note that Theorem 4.6 is also important for proving
regularity results for the Green function associated with the Lam\'e system.

(3) Theorem 4.6 should have important applications in the numerical treatment
of (4.5)-(4.6). For instance, it may be used to establish error estimates for
the so called collocation method.

(4) Another remark of interest is that for Lipschitz domains with a more
complicated underlying topology most boundary integral operators are only
Fredholm. In this case, Theorem 2.9 can be utilized. See \cite{30}.

(5) Finally, we want to stress that our techniques point to somewhat similar
results for the heat operator $\partial_t-\triangle$ in Lipschitz cylinders
as well. The corresponding $H^1$-atomic theory for the adjoint of the caloric
double layer potential operator has been worked out in \cite{3}.

\vskip 0.20in
\noindent --------------------------------------
\vskip 0.8in

\noindent Nigel Kalton

\noindent Department of Mathematics

\noindent University of Missouri-Columbia

\noindent Columbia, MO 65211

\vskip 0.15in

\noindent Marius Mitrea

\noindent Department of Mathematics

\noindent University of Missouri-Columbia

\noindent Columbia, MO 65211

\vskip 0.10in

\noindent and
\vskip 0.10in

\noindent The Institute of Mathematics

\noindent of the Romanian Academy,

\noindent P.O. Box 1-764

\noindent RO-70700 Bucharest, Romania

\enddocument